\documentclass[preprint,12pt]{elsarticle}
\usepackage{amssymb,amsmath,multirow}
\journal{ }

\begin{document}
%%%%%%%%%%%%%%%%
\begin{frontmatter}



\title{An Improved Analysis of Semidefinite Approximation
Bound for Nonconvex Nonhomogeneous Quadratic
Optimization with Ellipsoid Constraints}
\tnotetext[0]{This research was supported by Beijing Higher Education Young Elite Teacher Project (No.29201442), the fund of State Key Laboratory of Software Development Environment (No.SKLSDE-2013ZX-13), and by the Chinese NSF under the grant 11101261.}
\author[ad1]{Yong Hsia}%$^1$
\ead{dearyxia@gmail.com}
\author[ad1]{Shu Wang}%$^1$
\ead{wangshu.0130@163.com}
\author[ad2]{Zi Xu\corref{cor1}}
\ead{xuzi@shu.edu.cn}
\cortext[cor1]{Corresponding author.}
\address[ad1]{State Key Laboratory of Software Development
              Environment, LMIB of the Ministry of Education, School of
Mathematics and System Sciences, Beihang University, Beijing,
100191, P. R. China}
\address[ad2]{Department of Mathematics, Shanghai University,  Shanghai
200444,  P. R. China}

\begin{abstract}
We consider the problem of approximating nonconvex quadratic optimization with ellipsoid constraints (ECQP).
We show some SDP-based approximation bounds for special cases of (ECQP) can be improved by trivially applying the extened Pataki's procedure. The main result of this paper is to  give a new analysis on approximating
(ECQP) by the SDP relaxation, which greatly improves Tseng's result [SIAM Journal Optimization, 14, 268-283, 2003].  As an application, we strictly improve the approximation ratio for the assignment-polytope constrained quadratic program.
% Enter your abstract
\end{abstract}%

%\KEYWORDS{deterministic inventory theory; infinite linear programming duality;
%  existence of optimal policies; semi-Markov decision process; cyclic schedule}


\end{frontmatter} 
%\HISTORY{This paper was
%first submitted on April 12, 1922 and has been with the authors for
%83 years for 65 revisions.}

%\maketitle
%%%%%%%%%%%%%%%%%%%%%%%%%%%%%%%%%%%%%%%%%%%%%%%%%%%%%%%%%%%%%%%%%%%%%%

% Samples of sectioning (and labeling) in OPRE
% NOTE: (1) \section and \subsection do NOT end with a period
%       (2) \subsubsection and lower need end punctuation
%       (3) capitalization is as shown (title style).
%
%\section{Introduction.}\label{intro} %%1.
%\subsection{Duality and the Classical EOQ Problem.}\label{class-EOQ} %% 1.1.
%\subsection{Outline.}\label{outline1} %% 1.2.
%\subsubsection{Cyclic Schedules for the General Deterministic SMDP.}
%  \label{cyclic-schedules} %% 1.2.1
%\section{Problem Description.}\label{problemdescription} %% 2.

% Text of your paper here

\section{Introduction}\label{sect1}

In this paper, we consider the following nonconvex quadratic optimization problem with ellipsoid constraints:
\begin{align}
\min_{x\in \Bbb R^n}&\quad f(x)= x^TAx+2b^Tx \label{ECQP}\tag{ECQP}\\
{\rm s.t.}&\quad \|F^kx+g^k\|^2\le 1,~k=1,\ldots,m,\nonumber
\end{align}
where $A\in\Bbb {R}^{n\times n}$ symmetric, $F^k\in\Bbb {R}^{r^{k}\times n}$,  $b\in \Bbb {R}^n$, $g^k\in \Bbb {R}^{r^{k}}$, $r^{k}\geq 1$ and $\|\cdot\|$ denotes the Euclidean norm. Generally, this problem is NP-hard. To avoid trivial cases, we assume the Slater condition holds, i.e., the feasible region of \eqref{ECQP} has an interior point. With a proper transformation if necessary, we first make the following assumption.
\begin{ass}\label{as1}
The origin $0$ is in the interior of the feasible region of  \eqref{ECQP}, that is,
\[
\|g^k\|<1,~k=1,\ldots,m.
\]
\end{ass}
\eqref{ECQP} can be homogenized as
\begin{eqnarray}
 &\min_{x\in\Bbb R^{n+1}}& \sum_{i=1}^{n+1}\sum_{j=1}^{n+1}B_{ij}x_ix_j \\
&{\rm s.t.}& \sum_{i=1}^{n+1}\sum_{j=1}^{n+1}B_{ij}^kx_ix_j\le 0,~k=1,\ldots,m,\\
&& x_{n+1}=1,
\end{eqnarray}
where
\[
B=\left[\begin{array}{cc}A& b\\b^T& 0\end{array}\right],~
B^k=\left[\begin{array}{cc}(F^k)^TF^k& (F^k)^Tg^k\\(g^k)^TF^k& \|g^k\|^2-1\end{array}\right],~k=1,\ldots,m.
\]
By letting $X=xx^T$ and dropping the rank one constraint, the semidefinite programming relaxation of \eqref{ECQP} can be written as follows.
\begin{align}
\min&\quad  B\bullet X \nonumber\\
{\rm s.t.}&\quad B^k\bullet X\le 0,~k=1,\ldots,m,\label{SDP}\tag{SDP}\\
&\quad X_{n+1,n+1}=1,~X\succeq 0,~X\in\Bbb R^{(n+1)\times (n+1)}.\nonumber
\end{align}

In addition, we need to make the following assumption for \eqref{SDP} throughout this paper.
\begin{ass}\label{as2}
\eqref{SDP} has an optimal solution $X^*$.
\end{ass}
Let $v(\cdot)$
denote the optimal value of problem $(\cdot)$. Obviously, we have
\[
v({\rm SDP}) \le v({\rm ECQP}),
\]
and the equality holds if and only if rank$(X^*)=1$ with $X^*$ being an optimal solution of \eqref{SDP}.
Generally, the following theorem shows that \eqref{SDP} can also give a guaranteed-approximate solution for \eqref{ECQP}.
\begin{thm}[\cite{Tseng}]\label{ts}
Under Assumptions \ref{as1} and \ref{as2}, a feasible solution $x$ for \eqref{ECQP} can be generated in polynomial time satisfying
\begin{equation}
f(x)\le \frac{(1-\gamma)^2}{ (\sqrt{m}+\gamma)^2}\cdot v({\rm SDP}),\label{tsb}
\end{equation}
where $\gamma:=\max_{k=1,\ldots,m}\|g^k\|$.
\end{thm}

One special case of \eqref{ECQP} is that $b=0$, $g^k=0$ for $k=1,\ldots,m$ and $\sum_{k=1}^m(F^k)^TF^k$ is positive definite. It was shown in \cite{Ne} that in this case a feasible solution $x$ can be generated from \eqref{SDP} satisfying
\begin{equation}
f(x)\le\frac{1}{2\ln(2(m+1)\mu)}\cdot v({\rm SDP})\label{nv}
\end{equation}
with $\mu:=\min\{m+1,\max_{k=1,\ldots,m}{\rm rank}((F^k)^TF^k)\}$. In particular, when \eqref{ECQP} has a ball constraint, $\mu=\min\{m+1,n\}$. Also for this special case, Ye and Zhang (Corollary 2.6 in \cite{yz}) showed that a feasible solution $x$ satisfying
\[
f(x)\le\frac{1}{\min\{m-1,n\}}\cdot v({\rm SDP}),
\]
 can be found. For more detailed results related to this special case, we refer to the survey paper \cite{He}.
%Even in this case, (\ref{tsb}) improves (\ref{nv}) for $m\le 11$.

Another special case is that $A\preceq0, b=0$ but $\|g^k\|$ ($k=1,\ldots,m$) are allowed to be nonzero. It is shown in \cite{Ye}
that a feasible solution $\widetilde{x}$ can be randomly generated in this case such that
\begin{equation}
E(\widetilde{x}^TA\widetilde{x})\leq \frac{(1-\max_{k}{\|g^k\|})^2}{4\ln(4mn\cdot\max_{k}\{{\rm rank}\left( ({F^k})^TF^k\right )\}}\cdot v({\rm SDP}), \label{ye}
\end{equation}
where $E(\cdot)$ is the expectation function. To be mentioned, the $n$ in the denominator should be $n+1$ according to the proof in \cite{Ye}.

This paper is organized as follows. By directly applying the extended Pataki's procedure, i.e., the algorithm RED in \cite{Be}, we show in Section 2 that both (\ref{nv}) and (\ref{ye}) can be further improved.
Our main result is shown in Section 3. We propose a sharper analysis on the semidefinte approximation bound for \eqref{ECQP}. More detailedly, from an optimal solution of \eqref{SDP}, a feasible solution $x$ for \eqref{ECQP} can be generated, which satisfies that
\[
f(x)\le \frac{(1-\gamma)^2}{\left(\sqrt{\widetilde{r}}+\gamma\right)^2}\cdot v({\rm SDP}),
\]
where $\widetilde{r}=\min\bigg\{\left\lceil\frac{ \sqrt{8m+17}-3}{2}\right\rceil,n+1 \bigg\}$ and $\gamma$ is defined the same as in Theorem \ref{ts}. This bound improves the result shown in Theorem \ref{ts} in the order $m$, i.e., from $O(1/m)$ to $O(1/\sqrt{m})$.

Moreover, in Section 4, for a special case of \eqref{ECQP}, i.e., the assignment-polytope constrained QP problem (AQP), we show a strictly improved approximation bound compared to the result in \cite{FLYe}. Although, it is claimed in \cite{Ye} that this ratio can be improved from $1/{O(n^3)}$ to $1/{O(n^2\log(4n^4))}$,  the analysis technique therein only works for a special case of (AQP).
%Moveover, our  result can be extended to the well-known quadratic assignment problem.
At last, some conclusions are given.

{\bf Notations.} Throughout the paper, $A \succeq 0$ stands for the matrix $A$ is positive  semidefinite,
$A\bullet B=\sum_{i,j=1}^{n}a_{ij}b_{ij}$ is the inner product of two matrices $A, B$. Let
 $\Bbb R^n$ and $S^n_{+}$ be the $n$-dimensional vector space and $n\times n$ positive semidefinite symmetric
matrix space, respectively.  The notation ``$:=$ '' denotes ``define''.

\section{Improved Approximation Bound for Two Special Cases}

In this section, two special cases of \eqref{ECQP} are considered.
Before giving the main results, we first restate the following key theorem given in \cite{Be} and omit the proof.
\begin{thm}[\cite{Be}]\label{thm:1}
Let $r$ be a positive integer. Suppose that (SDP) is solvable and
\begin{equation}
m+1\le (r+2)(r+1)/2-1. \label{rela}
\end{equation}
Then (SDP) has a solution $X^*$ for which rank($X^*$) $\le r$.
\end{thm}

It can be easily verified that (\ref{rela}) is equivalent to
\begin{equation}
r\ge  \left\lceil\frac{ \sqrt{8m+17}-3}{2}\right\rceil:=r_{0}.\label{r}
\end{equation}

Moreover, an algorithm called ``algorithm RED'' is proposed in \cite{Be} to find such a solution with rank less than or equal to $r_0$. This algorithm can be regarded as an extension of Pataki's procedure [\cite{GP,G.P}].

{\bf Case I}: Let $g^k=0$ for $k=1,\ldots,m$, and assume $\sum_{k=1}^m(F^k)^TF^k$ is positive definite. In this case, by using \eqref{r}, we can improve the result given in \cite{Ne} to be as follows.

\begin{thm}
Let $X^*$ be an optimal solution of \eqref{SDP} with $rank (X^*) \leq r_0$, then a feasible solution $x$ can be generated from $X^*$, and we have
\[
f(x)\le\frac{1}{2\ln(2(m+1)\bar\mu)}\cdot v({\rm SDP}),\label{nb}
\]
where $\bar\mu:=\min\{r_{0}+1,\max_{k=1,\ldots,m}{\rm rank}((F^k)^TF^k)\}$ and $r_0$ is given in \eqref{r}.
\end{thm}

Since the proof of this theorem is almost the same as that in \cite{Ne} except that we use an optimal solution of \eqref{SDP} with the rank being less than or equal to $r_0$ by \eqref{r} instead of $m$, we omit the detail here.

{\bf Case II}: We assume $A\preceq0, b=0$. Similar to Case I, by using \eqref{r}, we can improve the approximation bound for the SDP relaxation that given in \cite{Ye}. The new result is shown in the following theorem and the proof is omitted too.

\begin{thm}\label{thm:2}
Let $X^*$ be an optimal solution of \eqref{SDP} with $rank (X^*) \leq r_0$, then a feasible solution $\tilde{x}$ can be generated from $X^*$, and the expectation of the objective satisfies that
\[
E (\widetilde{x}^TA\widetilde{x})\leq \frac{(1-\max_{k}{\|g^k\|})^2}{4\ln(4m\widetilde{r}\overline{r})}\cdot v({\rm SDP}),
\]
where $\widetilde{r}=\min\{r_{0},n+1\}$,$\overline{r}=\min\{r_{0}, \max_{k}\{{\rm rank}({F^k})^TF^k)\}\}$ and $r_{0}$ is given in \eqref{r}.
\end{thm}

\section{Improved Approximation Bound for General Case}

In this section, we consider \eqref{ECQP} in general case. We aim to analyze the approximation bound for \eqref{SDP}. Before giving the main result, we first introduce the following theorem proposed in \cite{sz}.
\begin{thm}[\cite{sz}]\label{thm:3}
Let $X$ be a positive semidefinite matrix of rank $r$. Then, $B\bullet X\le0$ if and only if there is a rank-one decomposition
\[
X=\sum_{i=1}^rw_iw_i^T
\]
such that $w_i^TBw_i\le 0$ for $i=1,\ldots,r$.
\end{thm}

Let $X^*$ be an optimal solution of (SDP) and $r$ be the rank of $X^*$. According to Theorem \ref{thm:1}, we can assume $r$ satisfies (\ref{r}).

Since $X^*_{n+1,n+1}=1$, it can be easily checked that $B^*\bullet X^*=0$ with
\[
B^*=\left[\begin{array}{cc}A& b\\b^T& -v({\rm SDP})\end{array}\right].
\]
It follows from Theorem \ref{thm:3} that there are vectors $w_i=(u_i^T,t_i)^T\in \Bbb R^{n}\times\Bbb R$, $i=1,\ldots,r$ such that
\[
X^*=\sum_{i=1}^rw_iw_i^T,~{\rm and}~w_i^TB^*w_i\le 0,~i=1,\ldots,r.
\]
Therefore, we obtain
\begin{eqnarray}
u_i^TAu_i+2t_ib^Tu_i &\le& v({\rm SDP})t_i^2,~i=1,\ldots,r,\label{10}\\
\sum_{i=1}^r \|F^ku_i+t_ig^k\|^2&=&B^k\bullet X^*+X^*_{n+1,n+1}\le1,~k=1,\ldots,m,\label{13}\\
\sum_{i=1}^rt_i^2&=&X^*_{n+1,n+1}=1.\label{12}
\end{eqnarray}
It follows from (\ref{13}) that
\begin{equation}
\|F^ku_i+t_ig^k\|^2/ t_i^2  \le   1/ t_i^2,~i=1,\ldots, r,~k=1,\ldots,m,\label{ineq}
\end{equation}
where $1/0:=+\infty$.

Then, according to (\ref{ineq}) and (\ref{12}),    we have
\begin{eqnarray*}
\min_{i=1,\ldots,r} \left\{ \max_{k=1,\ldots,m}\|F^ku_i+t_ig^k\|^2/t_i^2\right\} &\le&
\min_{i=1,\ldots,r} \left\{1/t_i^2\right\}\\
&\le& \max_{\sum_{i=1}^{r}t_i^2=1} \left\{\min_{i=1,\ldots,r} \left\{1/t_i^2\right\}\right\} \\
&=& \max_{\sum_{i=1}^{r}t_i^2=1,~z\le 1/t_i^2,~i=1,\ldots,r } \left\{z\right\}\\
&\le& \max_{\sum_{i=1}^{r}1/z\ge1} \left\{z\right\}\\
&=&r,
\end{eqnarray*}
where the last inequality actually can hold as an equality. Now, we have shown that there is an index $\overline{i}$ such that
\begin{equation}
\|F^ku_{\overline{i}}/t_{\overline{i}}+g^k\|\le \sqrt{r},~k=1,\ldots,m.\label{16}
\end{equation}
Define
\begin{eqnarray*}
\overline{x}&:=&\left\{\begin{array}{ll}u_{\overline{i}}/t_{\overline{i}},
& {\rm if}~b^Tu_{\overline{i}}/t_{\overline{i}}\le 0,\\
-u_{\overline{i}}/t_{\overline{i}},
& {\rm otherwise},\end{array}\right.\\
\overline{\tau}&:=&\max\{\tau\in[0,1]:~\|\tau F^k\overline{x}+g^k\|^2\le 1,~k=1,\ldots,m\}.
\end{eqnarray*}
Now, we are ready to present our main result shown in the following theorem, which improves Theorem \ref{ts} significantly. Though the remaining proof of Theorem \ref{main} is very similar to that of Theorem \ref{ts}, we state the theorem and provide the detail proof here for the sake of completeness.
\begin{thm}\label{main}
Under Assumptions \ref{as1} and \ref{as2}, the above construction gives a feasible solution $x=\overline{\tau}\overline{x}$ satisfying
\begin{equation}
f(x)\le \frac{(1-\gamma)^2}{\left(\sqrt{\widetilde{r}}+\gamma\right)^2}\cdot v({\rm SDP}),\label{new}
\end{equation}
where $\gamma=\max_{k=1,\ldots,m}\|g^k\|$ and $\widetilde{r}=\min\{r_{0},n+1\}$.
\end{thm}

{\bf {Proof.}}
To be mentioned, we only consider the case that $m$ is not very large, i.e., $m<n+1$, otherwise $\widetilde{r}=n+1$ and the approximation bound remains the same as that in Theorem \ref{ts}.
We first estimate $\overline{\tau}$. Fix any $k\in \{1,\ldots,m\}$. Then from \eqref{16}, we can get that
$\|F^k\overline{x}+g^k\|\leq \sqrt{r}$ if $b^Tu_{\overline{i}}/t_{\overline{i}}\le 0$. Otherwise, we
have
\[
\|F^k\overline{x}+g^k\|=\|-(F^ku_{\overline{i}}/t_{\overline{i}}+g^k)+2g^k\|\leq \sqrt{r}+2\|g^k\|.
\]
Therefore, for any $\tau\in[0,1]$, we obtain
\[
\|F^k(\tau\overline{x})+g^k\|=\|\tau(F^k\overline{x}+g^k)+(1-\tau)g^k\|\leq \tau(\sqrt{r}+2\|g^k\|)+(1-\tau)\|g^k\|.
\]
Whenever $\tau\leq (1-\|g^k\|)/(\sqrt{r}+\|g^k\|)$, it can be easily checked that $\|F^k(\tau\overline{x})+g^k\|\leq 1$ since
$\|g^k\|\leq1$. Thus,
\[
\overline{\tau}\geq\min_{k=1,\ldots,m}{\frac{1-\|g^k\|}{\sqrt{r}+\|g^k\|}}
=\frac{1-\max_{k=1,\ldots,m}\|g^k\|}{\sqrt{r}+\max_{k=1,\ldots,m}\|g^k\|},
\]
where the equality is due to the fact that $f(\gamma)=(1-\gamma)/(\sqrt{r}+\gamma)$ is a decreasing function for
$\gamma\in [0,1)$.\par
Since $\tau\in[0,1]$, we have $\tau\geq \tau^2$ and thus
\begin{align}
f(\tau\overline{x}) =&\tau^2\overline{x}^TA\overline{x}+\tau b^T\overline{x}\nonumber\\
 \leq &\tau^2\overline{x}^TA\overline{x}+\tau^2 b^T\overline{x}\nonumber\\
 \leq &\tau^2\overline{x}^TA\overline{x}+\tau^2 b^Tu_{\overline{i}}/t_{\overline{i}}\label{add1}\\
 =& \tau^2({u_{\overline{i}}}^TAu_{\overline{i}}+t_{\overline{i}}b^Tu_{\overline{i}})/{t_{\overline{i}}}^2\nonumber\\
 \leq &\tau^2v({\rm SDP})\label{add2},
\end{align}
where \eqref{add1} holds because $b^T\overline{x}\leq b^Tu_{\overline{i}}/t_{\overline{i}}, $ which is implied by the choice of $\overline{x}$,
and \eqref{add2} follows from (\ref {10}). By Assumptions \ref{as1}, 0 is a feasible solution to \eqref{ECQP} and hence $v({\rm SDP})\leq v({\rm ECQP})\leq f(0)=0$. Then the proof is completed if we set $\tau=\overline{\tau}$.
{\hfill \quad{$\Box$} }

Notice that Theorem \ref{main} remains the same as Theorem \ref{ts} when $m=1,2$ since
\[
m= \left\lceil\frac{ \sqrt{8m+17}-3}{2}\right\rceil {\rm ~when~} m=1,2.
\]
However, it strictly improves Theorem \ref{ts} when $m\ge 3$ since
\[
m> \left\lceil\frac{ \sqrt{8m+17}-3}{2}\right\rceil {\rm ~when~} m\ge 3.
\]
 For the special case that $g^k=0$ for $k=1,\ldots,m$ and there is a $k$ such that $(F^k)^TF^k$ is positive definite, our bound (\ref{new}) strictly improves (\ref{nv}) when $m\le 323$.

 \section{Application to the assignment-polytope constrained quadratic program}
In this section, we consider the following assignment-polytope constrained quadratic program:
\begin{align}
\min&\quad f(x)= x^TAx+2b^Tx \tag{ASQP}\\
{\rm s.t.}&\quad x\in F,\nonumber
\end{align}
where $F=\{x\in \Bbb R^{n^2}: \sum_{i=1}^n x_{i,j}=1,\sum_{j=1}^n x_{i,j}=1, i,j=1,\ldots,n, x_{i,j}\geq 0\}$.
Denote by $\overline{p}$ and $\underline{p}$ the maximal and minimal objective values $f(x)$ over $F$, respectively. Then,
an $\epsilon$-minimal solution  ($\epsilon \in [0,1]$) for (\rm ASQP) is defined as an $x\in F$ such that
\[
\frac{f(x)-\underline{p}}{\overline{p}-\underline{p}} \leq \epsilon.
\]
Fu et al.\cite{FLYe} showed that a $\left(1-\frac{1}{n^2(2n-2)}+\frac{1}{n^3(2n-2)}\right)$-minimal solution can be found in polynomial time.

%From the first two equality constraints, we can get
Since all the vectors satisfying the equality constraints in $F$  can be expressed as
\begin{equation}
x=\frac{1}{n}e+Ny,~ y\in \Bbb R^{(n-1)^2},\label{x-y}
\end{equation}
where $N\in \Bbb R^{n^2\times(n-1)^2}$ is the matrix basis of the null space for the equality constraints,
and $e\in \Bbb R^{n^2}$ is the vector of all ones,  the feasible region  $F$ in terms of $y$ becomes
\begin{equation}
-\frac{1}{n}\leq N_{i}y\leq 1-\frac{1}{n},i=1,\ldots,n^2,\label{ny}
\end{equation}
where $N_{i}$ is the $i$th row of $N$.

Now, we can reformulate (ASQP) as instances of (ECQP) in terms of $y$:
\begin{align}
\min& \quad h(y)= \left(\frac{1}{n}e+Ny\right)^TA\left(\frac{1}{n}e+Ny\right)+2b^T\left(\frac{1}{n}e+Ny\right) \tag{ASQP'}\\
{\rm s.t.}&\quad \left(2N_{i}y+\left(\frac{2}{n}-1\right)\right)^2\leq 1,~i=1,\ldots,n^2.\nonumber
\end{align}
As a corollary of (\ref{ye}), Ye \cite{Ye} gave an approximation algorithm, which generates a feasible point $\widetilde{y}$ such that \[
E(\widetilde{y}^TN^TAN\widetilde{y})\le\frac{1}{n^2\log(4n^4)}\cdot v{\rm (ASQP')},
\]
under the assumption that $h(y)$ is homogeneous and $N^TAN\preceq 0$. We notice that Ye's result is very special  since $h(y)$ is nonhomogeneous even when $f(x)$ is homogeneous.

Before applying Theorem \ref{main} to (ASQP'), we can easily see that
\[
\gamma=\max_{k=1,\ldots,m}\|g^k\| =1-\frac{2}{n}\ge0
\]
as $n\ge 2$.
%which can be substituted in the expression (\ref{new}).
Then, it follows from Theorem \ref{main} that  we can find a feasible solution $y$ such that
\begin{equation}
h(y)-h(0)\le g(n)\cdot \left(v({\rm SDP})-h(0)\right)\le g(n)\cdot \left(v{\rm (ASQP')}-h(0)\right),\label{new2}
\end{equation}
where
\[
g(n):=\frac{4}{n^2\left(\sqrt{\left\lceil\frac{ \sqrt{8n^2+17}-3}{2}\right\rceil}+1-\frac{2}{n}\right)^2}.
\]
Now, for (ASQP), we have
\begin{cor}
We can find a (1-g(n))-minimizer of (ASQP) in polynomial time.
%\begin{equation}
%f(x)\le \frac{4}{n^2\left(\sqrt{\left\lceil\frac{ \sqrt{8n^2+17}-3}{2}\right\rceil}+1-\frac{2}{n}\right)^2}\cdot v({\rm SDP}):=\frac{1}{g(n)}v({\rm SDP}),\label{new}
%\end{equation}
\end{cor}

 {\bf {Proof.}}
%where $g(n)$ is a function of variable $n$.
%Comparing the result with Ye\cite{Ye}, the new result is more general.
%In fact, before comparing the two bounds, we first change
%the form of our bound as follows:
We fist find a vector $y$ satisfying (\ref{new2}) and then generate $x$ according to (\ref{x-y}). Since
\[
f(x)=h(y),~f\left(\frac{1}{n}e\right)=h(0)\le \overline{p},~v{\rm (ASQP')}=v{\rm (ASQP)}=\underline{p},
\]
it follows from (\ref{new2}) that
\begin{eqnarray*}
f(x) &\le& f\left(\frac{1}{n}e\right)+g(n)\cdot \left(\underline{p}-f\left(\frac{1}{n}e\right)\right)\\
&=&(1-g(n))\cdot f\left(\frac{1}{n}e\right)+g(n)\cdot \underline{p}\\
&\le&(1-g(n))\cdot \overline{p}+g(n)\cdot \underline{p}
\end{eqnarray*}
Therefore, it holds that
\[
\frac{f(x)-\underline{p}}{\overline{p}-\underline{p}}
\leq \frac{(1-g(n))\cdot \overline{p}+g(n)\cdot \underline{p}  -\underline{p}}{\overline{p}-\underline{p}}= 1- g(n).
\]
The proof is complete. {\hfill \quad{$\Box$} }

%Thus we can easily see that the $1-\frac{1}{g(n)}$ is actually equivalent to $\epsilon$. Then we can find the new bound (\ref{new})
%improves the bound of Fu et al. for $m\leq10$ obviously. When $m>10$, the two bounds are almost equal.
Our new bound strictly improves that of Fu et al. \cite{FLYe} since
\[
1-\frac{1}{n^2(2n-2)}+\frac{1}{n^3(2n-2)}>1-\frac{1}{2}g(n)>1-g(n),
\]
which can be verified as follows by noting $n\ge 2$:
\begin{eqnarray*}
g(n)&=&\frac{4}{n^2\left(\sqrt{\left\lceil\frac{ \sqrt{8n^2+17}-3}{2}\right\rceil}+1-\frac{2}{n}\right)^2}\\
&>&\frac{4}{n^2\left(\sqrt{\left\lceil  \sqrt{2}n-0.7\right\rceil}+1-\frac{2}{n}\right)^2}~~({\rm since}~\sqrt{8n^2+17}<\sqrt{8}n +1.6 )\\
&\ge&\frac{4}{n^2\left(\sqrt{   \sqrt{2}n+0.3} +1-\frac{2}{n}\right)^2}\\
&>&\frac{4}{n^2\left(\sqrt[4]{2}\sqrt{n}+1.1-\frac{2}{n}\right)^2} ~~({\rm since}~\sqrt{   \sqrt{2}n+0.3}<\sqrt[4]{2}\sqrt{n}+0.1 )\\
&>&\frac{4}{n^2\left(\sqrt[4]{2}\sqrt{n}+1.1\right)^2}\\
&>&\frac{4}{n^2\left(0.7\sqrt{8}\sqrt{n}\right)^2} ~~({\rm since}~\sqrt[4]{2}\sqrt{n}+1.1<0.7\sqrt{8}\sqrt{n})\\
&>&\frac{1}{n^3}\\
&=&2\left(\frac{1}{n^2(2n-2)}-\frac{1}{n^3(2n-2)}\right).
\end{eqnarray*}

\section{Conclusion}
In this paper, we have proposed an improved analysis of the semidefinite approximation bound for nonconvex quadratic optimization problem with ellipsoid constraints. Two special cases are also discussed. As an application, we strictly improves the approximation bound for the assignment-polytope
constrained quadratic program. It is still need to be further investigated whether the new bounds are tight or not.

% References here (outcomment the appropriate case)

% CASE 1: BiBTeX used to constantly update the references
%   (while the paper is being written).
%\bibliographystyle{ormsv080} % outcomment this and next line in Case 1
%\bibliography{<your bib file(s)>} % if more than one, comma separated

% CASE 2: BiBTeX used to generate mypaper.bbl (to be further fine tuned)
%\input{mypaper.bbl} % outcomment this line in Case 2

%If you don't use BiBTex, you can manually itemize references as shown below.

%%%%%%%%%%%%%%%%%
\end{document}